\documentstyle[12pt]{article}
\newcommand{\const}{\mathop{\rm const}\limits}

\newcommand{\supp}{\mathop{\rm supp}\limits}

\begin{document}

\begin{center}

\vspace{3mm}

{\bf SHARP LYAPUNOV'S INEQUALITY FOR THE }\par

\vspace{4mm}

{\bf MEASURABLE SETS WITH INFINITE MEASURE,} \\

\vspace{4mm}

{\bf with generalization to the Grand Lebesgue spaces.}

\vspace{4mm}

 $ {\bf E.Ostrovsky^a, \ \ L.Sirota^b } $ \\

\vspace{4mm}

$ ^a $ Corresponding Author. Department of Mathematics and computer science, Bar-Ilan University, 84105, Ramat Gan, Israel.\\

E-mail:  \  eugostrovsky@list.ru\\

\vspace{3mm}

$ ^b $  Department of Mathematics and computer science. Bar-Ilan University,
84105, Ramat Gan, Israel.\\

E-mail: \ sirota3@bezeqint.net \\

\vspace{4mm}
                    {\sc Abstract.}\\

 \end{center}

 \vspace{3mm}

  We extend the classical Lyapunov inequality  on the  measurable space with infinite measure and on the
so-called Grand Lebesgue spaces (GLS).\par

 We find also the exact value for correspondent constant. \par

 Possible applications: Functional Analysis (for instance, interpolation of operators), Integral Equations,
Probability Theory and Statistics (tail estimations for random variables) etc.\par

 \vspace{4mm}

 \par

\vspace{3mm}


\vspace{4mm}

\section{Notations. Statement of problem.}

\vspace{4mm}

 Let $ (X = \{x\}, {\bf B}, \mu) $ be measurable space with  non-trivial sigma-finite measure $ \mu. $ We denote as ordinary
for arbitrary measurable numerical (or complex valued) function $ f: X \to R $  and for any measurable set $ A: A \in {\bf B} $
with finite  non-zero measure: $  0 < \mu(A) < \infty $

$$
|f|_{p,A} = \left[  \int_A |f(x)|^p \ d \mu(x) \right]^{1/p},\eqno(1.1)
$$

$$
|f|_{p} = \left[ \int_X |f(x)|^p \ d \mu(x) \right]^{1/p}, \ 1 \le p < \infty,
$$

$$
f \in L_p(A) \ \Leftrightarrow |f|_{p,A} < \infty, \  f \in L_p \ \Leftrightarrow |f|_{p} < \infty.
$$

 Note that $  |f|_{p} = |f|_{p,X}. $ \par

 The classical Lyapunov's inequality asserts  that if  $ \mu(X) = 1, \ 1 \le p \le q,  $  then \\
 $ |f|_p \le |f|_q.  $ Therefore under at the same restrictions on the variables $ p,q $ there holds the following inequality

$$
\frac{|f|_{p,A}}{ \mu^{1/p}(A)} \le  \frac{|f|_{q,A}}{\mu^{1/q}(A) }, \hspace{5mm} 0 < \mu(A) < \infty. \eqno(1.2)
$$

\vspace{4mm}

 {\bf  Our purpose in this short report is to generalize the Lyapunov's inequality (1.2)
 into the so-called Grand Lebesgue Spaces (GLS) instead the classical Lebesgue-Riesz spaces
 constructed over the measurable space with (in general case) infinite measure.   }\par

\vspace{4mm}

 We must  recall here briefly  the definition and some simple properties of the so-called Grand Lebesgue spaces;   more detail
investigation of these spaces see in \cite{Fiorenza3}, \cite{Iwaniec2}, \cite{Kozachenko1}, \cite{Liflyand1}, \cite{Ostrovsky1},
\cite{Ostrovsky2}; see also reference therein.\par

  Recently  appear the so-called Grand Lebesgue Spaces $ GLS = G(\psi) =G\psi =
 G(\psi; a,b), \ a,b = \const, a \ge 1, a < b \le \infty, $ spaces consisting
 on all  the measurable functions $ f: X \to R $ with finite norms

$$
   ||f||G(\psi) \stackrel{def}{=} \sup_{p \in (a,b)} \left[ |f|_p /\psi(p) \right]. \eqno(1.3)
$$

  Here $ \psi(\cdot) $ is some continuous positive on the {\it open} interval
$ (a,b) $ function such that

$$
     \inf_{p \in (a,b)} \psi(p) > 0, \ \psi(p) = \infty, \ p \notin (a,b).
$$

 We will denote
$$
 \supp (\psi) \stackrel{def}{=} (a,b) = \{p: \psi(p) < \infty, \}
$$

The set of all $ \psi $  functions with support $ \supp (\psi)= (a,b) $ will be
denoted by $ \Psi(a,b),  $ and and denote the set of all such a functions as $  \Psi: $

$$
\Psi \stackrel{def}{=} \cup_{1 \le a < b \le \infty} \Psi(a,b).
$$

  The Grand Lebesgue spaces (GLS) are rearrangement invariant, see \cite{Bennet1}, \cite{Krein1}, \cite{Kufner1}, \cite{Lieb1},
and  are used, for example, in the theory of probability \cite{Kozachenko1},  \cite{Ostrovsky1}, \cite{Ostrovsky2}; theory of Partial
Differential Equations \cite{Fiorenza3}, \cite{Iwaniec2};  functional analysis \cite{Fiorenza3}, \cite{Iwaniec2},  \cite{Liflyand1},
 \cite{Ostrovsky2}; theory of Fourier series, theory of martingales, mathematical statistics, theory of approximation etc.\par

\vspace{3mm}

 The so-called {\it fundamental function} $ \phi(\delta) = \phi(G\psi(a,b), \delta), \ \delta \in (0, \infty)  $
  of the $ G\psi(a,b) $ space  may be calculated by the formula

$$
\phi(G\psi(a,b), \delta) = \sup_{p \in (a,b)} \left[ \frac{\delta^{1/p}}{\psi(p)}  \right].
$$

 This notion play a very important role in the theory of interpolation of operators, Fourier series
\cite{Bennet1},   theory of random variables, in particular, theory of Central Limit Theorem,
in Banach spaces \cite{Ostrovsky10} etc. \par
 The detail investigation of fundamental  functions for Grand Lebesgue Spaces with consideration of some examples
see in the articles \cite{Liflyand1}, \cite{Ostrovsky11}. \par

 Let us introduce also the  so-called {\it localized} GLS norm. Indeed,  we define for $ \psi \in \Psi(a,b), \
 A: A \in {\bf B}, \ \mu(A) \in (0, \infty), $ and measurable function $  f: A \to R $

 $$
 ||f||_{\psi,A} := \sup_{p \in (a,b) } \left[ \frac{|f|_{p,A}}{\psi(p)} \right] =
 \sup_{p \in (a,b) } \left[ \frac{[\int_A |f(x)|^p \mu(dx) ]^{1/p}}{\psi(p)} \right]. \eqno(1.4)
 $$

\vspace{4mm}

\section{Main result.}

\vspace{4mm}

 Let $ \psi, \nu $ be two function from the set $  \Psi. $

\vspace{3mm}

{\bf Definition  2.1.} The following functional $ R(G\psi, G\nu) = $

$$
R(\psi, \nu) \stackrel{def}{=} \sup_{ 0 \ne f \in G\psi \cap G\nu}
\sup_{A: 0 < \mu(A) < \infty} \left[ \frac{||f||_{\psi,A}}{\phi(G(\psi), \mu(A))}: \frac{||f||_{\nu,A}}{\phi(G(\nu), \mu(A))} \right] \eqno(2.1)
$$
will named  {\it double ratio}  between the spaces $ G\psi $ and $ G\nu. $ \par

\vspace{3mm}

{\bf Theorem 2.1.} {\it Let $ \psi, \ \nu $ be two function such that $ \psi \in \Psi(a_1,b_1), \ \nu \in \Psi(a_2, b_2). $
Suppose $ b_1 < a_2; $ the opposite case is trivial for us.  Our statement: }

$$
R(G\psi, G\nu) =  1. \eqno(2.2)
$$

\vspace{3mm}

{\bf Proof.} \\

 {\bf First step.} Let $ f \in G\psi \cap G\nu, f \ne 0, \ p \in (a_1,b_1), \ g \in (a_2,b_2), \ 0 < \mu(A) < \infty.  $
We rewrite the Lyapunov's inequality (2.1) as follows:

$$
|f|_{p,A}  \le  \frac{|f|_{q,A}}{\mu^{1/q}(A) } \cdot \mu^{1/p}(A),
$$
or after dividing over $  \psi(p) $

$$
\frac{|f|_{p,A}}{ \psi(p) }  \le  \frac{|f|_{q,A}}{\mu^{1/q}(A) } \cdot \frac{ \mu^{1/p}(A)}{ \psi(p) }. \eqno(2.3)
$$

 We have taking the supremum over $ p, \ p \in (a_1, b_1) $  from both the sides of inequality (2.3)
 using the direct definition of the Grand Lebesgue norm and also the definition of the fundamental function

$$
||f||_{\psi,A} \le \frac{|f|_{q,A}}{\mu^{1/q}(A) } \cdot \phi(G\psi,\mu(A)). \eqno(2.4)
$$

\vspace{3mm}

{\bf Second step.}  Further, we use the simple estimate $ |f|_{q,A} \le ||f||_{\nu,A}\cdot \nu(q), \ q \in (a_2,b_2), $
in the inequality (2.4):

$$
\frac{||f||_{\psi,A}}{ \phi(G\psi,\mu(A)) } \le ||f||_{\nu,A} \cdot \frac{\nu(q)}{\mu^{1/q}(A)}  =
$$

$$
||f||_{\nu,A} \cdot \left[ \frac{\mu^{1/q}(A)}{ \nu(q)  }  \right]^{-1},
$$
and after taking infinum over $  q: $

$$
\frac{||f||_{\psi,A}}{ \phi(G\psi,\mu(A)) } \le ||f||_{\nu,A} \cdot \inf_{q \in (a_2,b_2)} \left[ \frac{\nu(q)}{\mu^{1/q}(A)} \right]  =
$$

$$
||f||_{\nu,A} \cdot \left[ \sup_{q \in (a_2,b_2)} \frac{\mu^{1/q}(A)}{\nu(q)} \right]^{-1} = \frac{||f||_{\nu,A}}{\phi(G\nu,\mu(A))},
$$
hence $  R(G\psi, G\nu) \le 1. $ \\

\vspace{3mm}

{\bf Lower bound.} It is easy to verify that the lower bound for the expression for the functional $ R(G\psi, G\nu) $
is attained, for example, if $  f_0(x) = 1, $ so that $  R(G\psi, G\nu)  \ge 1. $ \par

In detail,

$$
\frac{||f_0||_{\psi,A}}{\phi(G\psi,\mu(A))} =
\frac{\sup_{p \in (a_1,a_2)} \left\{ \left[ \int_A 1 \cdot \mu(dx)  \right]^{1/p}/\psi(p) \right\} }{\phi(G\psi,\mu(A))} =
$$

$$
\frac{\sup_{p \in (a_1,a_2)} \left\{ \left[ \mu(A)\right]^{1/p}/\psi(p) \right\} }{\phi(G\psi, \mu(A))} =
\frac{\phi(G\psi, \mu(A))}{\phi(G\psi, \mu(A))}= 1
$$
and analogously $ ||f_0||_{\nu,A} /\phi(G\nu,A) = 1.$ \par

\vspace{3mm}

 This completes the proof of theorem 2.1.\\

\vspace{4mm}

{\bf Remark 2.1.} The definition (2.1) may be easily extended on  arbitrary pair of rearrangement invariant (r.i.) spaces
$ F_1  $ and $ F_2 $ with norms correspondingly $  ||f||F_1, \ ||f||F_2,  $
constructed over source triple $  (X, {\bf B,} \mu). $ Namely, denote

$$
||f||_{F_1,A} = ||f \cdot I_A||F_1, \hspace{6mm}  ||f||_{F_2,A} = ||f \cdot I_A||F_2,
$$
where as usually $ I_A = I_A(x), \ x \in X $ is indicator function  generated by the measurable set $  A. $ Then by definition of
the double ratio  between the spaces $ F_1 $ and $ F_2 $  is the following functional

$$
R(F_1, F_2) \stackrel{def}{=} \sup_{ 0 \ne f \in F_1 \cap F_2}
\sup_{A: 0 < \mu(A) < \infty} \left[ \frac{||f||_{F_1,A}}{\phi(F_1, \mu(A))}: \frac{||f||_{F_2,A}}{\phi(F_2, \mu(A))} \right].
$$

 Obviously,  $ R(F_1, F_2) \ge 1 $ for any pair of r.i. spaces $ F_1, \ F_2. $ \par

\vspace{3mm}

{\bf Remark 2.2.} If we introduce the {\it discontinuous} function

$$
\psi_r(p) = 1, \ p = r; \psi_r(p) = \infty, \ p \ne r, \ p,r \in (a,b)
$$
and define formally  $ C/\infty = 0, \ C = \const \in R^1, $ then  the norm
in the space $ G(\psi_r) $ coincides with the $ L_r $ norm:

$$
||f||G(\psi_r) = |f|_r.
$$

 Thus, the inequality (1.2) is particular, more exactly, extremal case of the assertion of theorem 2.1.\\

\vspace{4mm}

\section{Concluding remarks.}

 \vspace{4mm}

 Our considerations are very similar to ones in the article \cite{Ostrovsky2}, devoted to the generalization of
Nilol'skii inequality, and may be used perhaps in   the theory of operators acting in the Lorentz spaces, see, e.g.
\cite{Astashkin1}, \cite{Osekowski1}, \cite{Osekowski2}, \cite{Kufner2}, \cite{Okikiolu1}.\\

\end{document}